\documentclass[smallextended,referee,envcountsect,]{svjour3}
\smartqed
\usepackage{graphicx}

\usepackage{float}
\usepackage{cite}
\usepackage{amscd} 
\usepackage{verbatim}
\usepackage{amsxtra}
\usepackage{csquotes}
\usepackage{mathtools}
\usepackage{scrextend}
\usepackage{amsmath,etoolbox}
\usepackage{amssymb}
\usepackage{paralist}
\usepackage{graphics} 
\usepackage{epsfig} 
\usepackage{subcaption}
\usepackage{graphicx}  
\usepackage{epstopdf}
\usepackage[colorlinks=true]{hyperref}
\usepackage{centernot}
\usepackage{extarrows}
\usepackage{mathrsfs}
\hypersetup{urlcolor=blue, citecolor=red}
\usepackage{microtype}
\usepackage[mathscr]{euscript}
\usepackage[shortlabels]{enumitem}
\usepackage{lipsum,calc}
\usepackage[linesnumbered,ruled,vlined]{algorithm2e}
\SetKwInput{KwInput}{Input}                
\SetKwInput{KwOutput}{Output}
\SetKwInput{KwWrite}{Write}
\SetKw{KwBy}{by}
\makeatletter
\newenvironment{sqcases}{%
  \matrix@check\sqcases\env@sqcases
}{%
  \endarray\right.%
}
\def\env@sqcases{%
  \let\@ifnextchar\new@ifnextchar
  \left[
  \def\arraystretch{1.2}%
  \array{@{}l@{\quad}l@{}}%
}
\makeatother

\makeatletter
\newcommand{\raisemath}[1]{\mathpalette{\raisem@th{#1}}}
\newcommand{\raisem@th}[3]{\raisebox{#1}{$#2#3$}}
\makeatother

\def\gi{\gamma_i}

\def\f{\longrightarrow}
\def\im{\Longrightarrow}

\def\N{\mathbb{N}}

\def\e{\varepsilon}
\def\l{\lambda}

\def\x{\bar{x}}

\def\u{\bar{u}}
\def\vv{\bar{v}}
\def\a{\alpha}

\def\<{\langle}
\def\>{\rangle}

\def\e{\varepsilon}
\def\R{\mathbb{R}}

\def\inte{\textnormal{int}\,}
\def\clo{\textnormal{cl}\,}

\def\bdry{\textnormal{bdry}\,}

\def\conv{\textnormal{conv}\,}

\def\gk{\gamma_k}
\def\gkn{\gamma_{k_n}}

\def\CO{\mathcal{C}}
\def\c{\mathsf{c}}
\def\-{\textnormal{-}}
\def\I{\mathcal{I}}
\newcommand*{\tran}{^{\mkern-1.5mu\mathsf{T}}}

\def\sp{\hspace{0.015cm}}
\def\bp{\hspace{-0.08cm}}
\def\bbp{\hspace{-0.04cm}}

\def\c{\bbp c}

\def\U{\mathscr U}

\def\0gk{\scaleto{0,\gk}{3.5pt}}

\binoppenalty=\maxdimen 
\relpenalty=\maxdimen 

\begin{document}

\title{Numerical Method for a Controlled Sweeping Process with Nonsmooth Sweeping Set}


\author{Chadi Nour \and Vera Zeidan }

\institute{Chadi Nour \at
             Department of Computer Science and Mathematics \\
             Lebanese American University\\
             Byblos Campus, P.O. Box 36\\
              Byblos, Lebanon\\
              cnour@lau.edu.lb 
           \and
              Vera Zeidan,  Corresponding author  \at
              Department of Mathematics\\
              Michigan State University\\
              East Lansing, MI 48824-1027, USA\\
              zeidan@msu.edu
}

\date{Received: date / Accepted: date}

\maketitle

\begin{abstract} The numerical method developed in \cite{verachadinum}  for optimal control problems involving sweeping processes with {\it smooth} sweeping set $C$ is generalized to the case where $C$ is {\it nonsmooth}, namely, $C$ is the intersection of a finite number of sublevel sets of smooth functions. The novelty  of this extension resides in producing  for the general setting a different  approach,  since the one used for  the {\it smooth}  sweeping sets is not applicable here. 

\end{abstract}


\section{Introduction} 
{\it Sweeping processes} refer to a specific category of differential inclusions that incorporates the normal cone to a set called {\it sweeping set}. This distinctive feature leads  to  differential inclusions that are {\it unbounded} and {\it discontinuous}.  The initial appearance of such a  model  goes back to  the papers \cite{moreau1,moreau2,moreau3} by J.J. Moreau in which he introduced this system as  a framework for investigating the dynamics of plasticity and friction. Subsequently, various adaptations of this model have surfaced in a wide range of applications, including but not limited to engineering, mechanics, crowd motion problems, and economics, etc (see \cite{outrata} and its references). 

Over the past few years, extensive research has been conducted on optimal control problems over various versions of sweeping processes with particular focus on establishing  the existence of optimal solutions and on  deriving  necessary optimality conditions, see e.g.,\cite{brokate,ccmn,ccmnbis,cmo0,cmo,cmo2,chhm2,chhm,cmn0,henrion,palladino,pinho,pinhoEr,pinholast,VCpaper,verachadisvaa,verachadijune,verachadijca,verachadi}. However,  numerical methods for such problems  are quite limited in the literature, with a few notable exceptions given in   \cite{outrata,pinhonum,verachadinum}. 

In this paper we are interested in constructing a numerical algorithm to solve a {\it fixed} time Mayer problem $(P)$ in which the dynamic is a {\it controlled sweeping process} $(D)$, the sweeping set $C$ is the intersection of the zero-sublevel sets of a finite sequence of functions $\psi_i\colon\R^n\f\R$, $i=1,\dots,r$, and the initial state is a fixed point $x_0$ in $C$. This problem $(P)$ was successfully treated {\it numerically} in \cite{pinhonum} for the {\it special} case: the initial state $x_0$ lies in the interior of $C$, and $C$ is convex and of class $\CO^2$, that is, $r=1$ and $\psi:=\psi_1$ is convex and $\CO^2$. The main idea used in \cite{pinhonum}, which differs from that used in \cite{outrata}, consists of approximating  ($D$) by the system ($D_{\gk}$), obtained by replacing in ($D$) the normal cone $N_C(x(t))$  by the penalty term $\gk e^{\gk\psi(x(t))} \nabla\psi(x(t))$. Then, the so obtained standard optimal control problem is solved numerically over {\it piecewise constant} controls. This numerical method  is  generalized in \cite{verachadinum}  to allow the initial state $x_0$ to be {\it any} point in $C$, including its boundary, and to {\it nonconvex} and $\CO^{1,1}$ sweeping sets $C$. However, the {\it smoothness} of $C$ remains an essential assumption in  \cite{verachadinum} which naturally excludes a large class of {\it nonsmooth} sweeping sets arising from applications. 

The goal of this paper is to expand the domain of applicability of the numerical method in \cite{verachadinum} to a general form of $(P)$, namely, for $r>1$. In this  case, $C$ is {\it not} necessarily smooth for being the intersection of the zero-sublevel sets of a finite sequence of $\CO^{1,1}$-functions $\psi_i$, $i=1,\dots,r$.  While transitioning from $r=1$ to $r>1$ might be initially perceived as a minor generalization, a close examination shows that  this is not the case,  since this transition  actually necessitates a {\it major} overhaul of the  approach used in \cite{verachadinum}. This is due to the fact that    when $r>1$, the set $C$ is the zero-sublevel of the function $\psi:=\max\{\psi_i : i=1,\dots,r\},$ which is {\it only} guaranteed to be {\it Lipschitz}, and hence,  it renders the exponential penalization technique of \cite{verachadinum}  inapplicable. To circumvent this major obstacle, we approximate the nonsmooth max-function $\psi$ by a well constructed sequence $\psi_{\gk}$ of $\CO^{1,1}$ functions, and hence, we use in the definition of $(D_{\gk})$ the exponential penalization technique for $\psi_{\gk}$.  It turns out   that the so-obtained  $(D_{\gk})$ is {\it equivalent} to an approximating control system having $r$-penalty terms that  involve $\nabla\psi_i$, $i=1,\dots,r$. 
However, having  solved the nonsmoothness issue with $\psi$  by means of $\psi_{\gk}$, we now encounter a new hurdle caused by  the generalized Hessian of the sequence  $(\psi_{\gk})_k$  {\it not} being uniformly bounded. This issue requires  new ideas that will be revealed when  establishing, parallel to  \cite{verachadinum}, the theoretical results  needed for the development of our numerical method. 

The layout of the paper is as follows. In the next section, we present our basic notations and definitions, and  we state our optimal control problem $(P)$ over a sweeping process.  In Section \ref{hypo}, we list our  hypotheses, and provide some preparatory results. In Section \ref{key}, we establish three theoretical results, namely, Propositions \ref{propo1}, \ref{propo2} and \ref{propo3}, that form the backbone of the main result obtained in Section \ref{algo} as  Theorem \ref{th1}. This theorem confirms that optimal trajectories of  well-constructed approximating problems  converge to an optimal trajectory for the original problem $(P)$.  This result leads to designing  our proposed  numerical algorithm in Section \ref{algo}. The effectiveness of this algorithm is tested on a numerical example elaborated in Section \ref{example}. More precisely, using our algorithm we compute a numerical optimal trajectory for our example and we show that it is actually a good approximation of an exact optimal trajectory produced by means of  the maximum principle established in \cite{verachadijca}. The last section of the paper contains some concluding remarks.

\section{Basic Notations and Definitions, and Statement of $(P)$}

\subsection{Basic Notations and Definitions}

We denote by $\|\cdot\|$ and $\<\cdot,\cdot\>$, the Euclidean norm and the usual inner product, respectively. The open and the closed unit balls are respectively denoted by $B$ and  $\bar{B}$. For $x\in\R^n$ and $\rho\geq0$, the open and the closed balls of radius $\rho$ centered at $x$ are respectively written as $B_{\rho}(x)$ and $\bar{B}_{\rho}(x)$. For a set $S \subset \R^n$, $\inte S$, $\bdry S$, $\clo S$, $\conv S$, and  $S^c$ designate  the interior, the boundary, the closure, the convex hull, and the complement of $S$, respectively. The Lebesgue space of essentially bounded measurable functions $f\colon [a,b]\f\R^n$ is denoted by $L^{\infty}([a,b];\R^n)$. For the sets of absolutely continuous functions and  of bounded variations  functions  $f\colon [a,b]\f\R^n$ we use, respectively,   $AC([a,b];\R^n)$ and $BV([a,b];\R^n)$.  A function $f\colon\R^n\f\R$ is $\CO^{1,1}$ if it is Fr\'echet differentiable with locally Lipschitz derivative. A function $\pi\colon A\f B$ is bi-Lipschitz if it is a Lipschitz bijection onto $B$, and its inverse is also Lipschitz. 

Now we present some notations and definitions from {\it Nonsmooth Analysis} and {\it Geometry}. For standard references, see, e.g.,  the monographs \cite{brudnyi,clarkeold,clsw,delfour,mordubook,rockwet,ThibaultBook}. Let $S$ be a nonempty and closed subset of $\R^n$, and let $s\in S$. The {\it proximal}, the {\it Mordukhovich} ({also known as \it limiting}), and the {\it Clarke normal} cones to $S$ at $s$ are denoted by $N^P_S(s)$, $N_S^L(s)$, and $N_S(s)$, respectively. For the {\it Clarke tangent} cone to $S$ at $s$, we use $T_S(s)$. The  set $S$ is said to be {\it epi-Lipschitz} if for all $s\in S$, the Clarke normal cone of $S$ at $s$ is {\it pointed}, that is, $N_S(s)\cap-N_S(s)=\{0\}$. For $\rho > 0$, the set $S$ is said to be $\rho$-\textit{prox-regular} if for all $s\in S$ and for all unit vector $\zeta\in N_S^P(s)$, we have $\langle \zeta,x-s\rangle \leq \frac{1}{2\rho} \|x-s\|^2$ for all $x\in S$. This latter inequality is known as the {\it proximal normal inequality}. Finally, we say that $S$ is {\it quasiconvex} if there exists $c\geq 0$ such that any two points $s_1,\sp s_2$ in $S$ can be joined by a polygonal line $\gamma$ in $S$ satisfying $l(\gamma) \leq c\sp \|s_1-s_2\|,$ where $l(\gamma)$ denotes the length of $\gamma$. 

\subsection{Statement of $(P)$}

 This paper focuses on developing a numerical algorithm to solve the following fixed time Mayer problem $$ \begin{array}{l} (P)\colon\; \hbox{Minimize}\;g(x(T))\\ \hspace{0.9cm} \hbox{over}\;(x,u)\;\hbox{such that}\;u\in \mathscr{U},\; x\in AC([0,T];\R^n),\, \hbox{and} \\[2pt]  \hspace{0.9cm} (D)  \begin{cases}\dot{x}(t)\in f(x(t),u(t))-N_C(x(t)), \;\;\hbox{a.e.}\;t\in[0,T]\\ x(0)=x_0, \end{cases}
    \end{array}$$
     where $T>0$ is fixed,  $g\colon\R^n\f\R$, $f\colon\R^n\times \R^m\f\R^n$,  $C$  is the intersection of the zero-sublevel sets of a finite sequence of functions $\psi_i\colon\R^n\f\R$, $i=1,\dots,r$, $N_C(\cdot)$ stands for the Clarke normal cone to $C$, $x_0\in C$ is fixed, and, for a given nonempty set $U\subset \R^m$, the set of control functions $\mathscr{U}$ is defined as
\begin{equation*}
\mathscr{U}:=\left\{ u\colon[0,T]\f\R^m \; \;  \hbox{is  measurable and}\;\; u(t) \in U \;\;\hbox{a.e.}\; t\in [0,T]\right\}. 
\end{equation*}
A pair $(x,u)$ is {\it admissible} for $(P)$ if $x\colon[0,T]\f\R^n$ is absolutely continuous, $u\in \mathscr{U}$, and $(x,u)$  satisfies the {\it controlled\,} {\it sweeping process} $(D)$ called the {\it dynamic} of $(P)$. An admissible pair $(\x,\u)$ for $(P)$ is said to be an {\it optimal solution} if $g(\x(T))\leq g(x(T))$ for all pairs $(x,u)$ admissible for $(P)$. In that case, $\x$  is called an {\it optimal trajectory} of $(P)$.

\section{Hypotheses and Preparatory Results}\label{hypo}

\subsection{Hypotheses} 
We assume throughout this paper that the data of $(P)$ satisfy the following hypotheses: 
\begin{enumerate}[label=\textbf{H\arabic*}:, leftmargin=\widthof{[H4]}+\labelsep]
\item $f$ is continuous on $C\times U$; and there exists $M>0$ such that $x\mapsto f (x, u)$ is $M$-Lipschitz for all $u\in U$; and $\|f(x,u)\| \leq M$ for all $(x,u)\in C \times U.$
\item $f(x,U)$ is convex for all $x\in\R^n$, and $U$ is compact.
\item $C\not=\emptyset$ is given by 
\begin{equation*} \label{Cdef} C:=\bigcap_{i=1}^{r}C_i,\;\hbox{where}\;C_i:=\{x\in\R^n : \psi_i(x)\leq0\}\;\hbox{and}\end{equation*} 
$(\psi_i)_{1\leq i\leq r}$ is a family of $\CO^{1,1}$ functions $\psi_i\colon\R^n\f\R$. Moreover, for $i=1,\dots,r$, $C_i$ is compact, with $\bdry C_i$ connected for $n>1$ and $C_i$ is convex for $n=1$, and  there is a constant $\eta>0$ such that \begin{equation*}\label{H3ineq} \left\|\sum_{i\in\I^0_c}\lambda_i\nabla\psi_i(c)\right\|>2\eta, \;\;  \forall c\in \{x\in \R^n: \I^0_x\not=\emptyset\},\end{equation*} where $\I^0_x:=\{i\in\{1,\dots,r\} : \psi_i(x)=0\}$ and $(\l_i)_{i\in\I^0_c}$ is any sequence of nonnegative numbers satisfying $\sum_{i\in\I^0_c}\l_i=1$.
\item $g$ is $L_g$-Lipschitz on $C$.
\end{enumerate}

We denote by $\bar{M}_\psi$ a {\it common upper bound} over $C$ of the finite sequence $(\|\nabla\psi_i(\cdot)\|)_{i=1}^{r}$ such that $\bar{M}_\psi\geq2\eta$, and by $2M_\psi$ a {\it common Lipschitz constant}  of the finite family $\{\nabla \psi_{i}\}_{i=1}^{r}$  over the compact set  $\conv C$. We also denote by:
\begin{enumerate}[$\bullet$]
\item $\psi\colon\R^n\f\R$ the function defined by
 \begin{equation*}\label{psidef} \psi(x):=\max\{\psi_i(x) : i=1,\dots,r\},\;\;\forall x\in\R^n. \end{equation*}
Clearly we have that $C=\{x\in\R^n : \psi(x)\leq 0\}.$
\item $(\gk)_{k}$ a sequence satisfying $\gk>\frac{2M}{\eta}$ for all $k\in\N$, with $\gk \f\infty\;\hbox{as}\;k \f\infty$.
\item $(\a_k)_k$ and $(\sigma_k)_k$ the two real sequences defined by\begin{equation}\label{sigmadef} \a_{k}:=\frac{\ln \left(\frac{\eta\gk}{2{M}}\right)}{\gk}\;\,\hbox{and}\;\,\sigma_k:=\frac{r{M}_{\psi}}{2\eta^2}\bigg(\frac{\ln(r)}{\gk}+\a_k\bigg),\;\;\;k\in\N. \end{equation}
For $(\a_{k})_k$, we have  $\gk e^{-\a_{k}\gk}=\frac{2{M}}{\eta},$ $\a_k>0$ for all $k\in\N$, and $\a_k\searrow 0$. For $(\sigma_k)_k$, we have  $\sigma_k>0$ for all $k\in\N$, and $\sigma_k\searrow 0$.
\item $(\psi_{\gk})_k$ the sequence of functions $\psi_{\gk}\colon\R^n\f\R$ defined by \begin{equation}\label{psigkdef} \psi_{\gk}(x):=\frac{1}{\gk}\ln\Bigg(\sum_{i=1}^{r}e^{\gk \psi_i(x)}\Bigg),\;\;\forall x\in\R^n. \end{equation}
Clearly we have that \begin{equation}\label{gradpsigk}
\nabla\psi_{\gk}(x)= \frac{\sum\limits_{i=1}^r e^{\gk\psi_i(x)} \nabla\psi_i(x)}{\sum\limits_{i=1}^re^{\gk\psi_i(x)}},\;\;\forall x\in\R^n.
\end{equation}
\item $(C_i(k))_k$ the sequence defined, for $i=1,\dots,r$, by 
 \begin{equation*}\label{Ci(k)def} C_i(k):=\{x\in\R^n : \psi_i(x)\leq -\a_k\}.\end{equation*}
\item $(C^{\gk})_k$ and $(C^{\gk}(k))_k$ the two sequences defined by  \begin{equation}\label{Cgkdef} C^{\gk}:=\{x\in\R^n : \psi_{\gk}(x)\leq 0\}=\bigg\{x\in\R^n : \sum_{i=1}^{r}e^{\gk \psi_i(x)}\leq 1\bigg\},\;\hbox{and}\end{equation}
 \begin{equation}\hspace{-0.4cm}\label{Cgkkdef} C^{\gk}(k):=\{x\in\R^n : \psi_{\gk}(x)\leq -\a_k\}=\bigg\{x\in\R^n : \sum_{i=1}^{r}e^{\gk \psi_i(x)}\leq \frac{2{M}}{\eta\gk}\bigg\}.\end{equation}
One can easily see that if $r=1$, then $\psi_{\gk}$, $C^{\gk}$, and $C^{\gk}(k)$ coincide with $\psi$, $C$, and $C(k):=\{x\in\R^n : \psi(x)\leq -\a_k\}$, respectively. 
\end{enumerate}

\subsection{Preparatory Results} This subsection consists of preparatory results that are fundamental for the rest of the paper. Note that some of these results are extracted from the papers  \cite{VCpaper,verachadijca,verachadi},  and hence, their proofs are omitted here. 

We begin with the following proposition which gives important properties of the set $C$.

\begin{proposition}[{\cite[Proposition 4.1]{verachadijca}}] \label{prop1} The set $C$ is $\frac{\eta}{M_\psi}$-prox-regular, epi-Lipschitz with $C=\clo(\inte C)$, and, for all $x\in\bdry C$ we have 
\begin{eqnarray*} \label{proxformulaprox} \nonumber N_C(x)=N^P_C(x)=N^L_C(x)&=&\bigcup_{\l\geq 0}\l\partial\psi(x)\\&=&\Bigg\{\sum_{i\in\I_x^0}\lambda_i\nabla\psi_i(x) : \l_i\geq 0\Bigg\}\not=\{0\}.\end{eqnarray*} 
Moreover, we have
$$\inte C=\bigcap_{i=1}^{r}\inte C_i=\bigcap_{i=1}^{r}\{x\in\R^n : \psi_i(x)<0\}\not=\emptyset,\;\;\hbox{and}$$ $$\bdry C=C\cap \bigg(\bigcup_{i=1}^{r} \bdry C_i\bigg)\not=\emptyset.$$
\end{proposition}

For the sequence of functions $(\psi_{\gk})_k$, we have the following.
\begin{proposition}[{\cite[Proposition 4.4]{verachadijca}}] \label{proppsigk} The following assertions hold\sp$:$
\begin{enumerate}[$(i)$, leftmargin=\widthof{(iii)}+\labelsep]
\item The sequence $(\psi_{\gk})_k\in \CO^{1,1}(\R^n;\R)$, is monotonically nonincreasing in  $k$, and converges uniformly to $\psi$. Moreover, for all $k\in\N$ and  for $i=1,\dots,r$, we have that
\begin{equation} \label{pisgkineq}\psi_i(x)\leq \psi(x)\leq \psi_{\gk}(x)\leq \psi(x)+\frac{\ln(r)}{\gk},\;\;\forall x\in\R^n,\;and
\end{equation}
\begin{equation*}\label{gradbdd}
\|\nabla\psi_{\gk}(x)\|\le \bar{M}_\psi, \;\; \forall x\in C.
\end{equation*}
\item  There exist $k_1\in\N$ and $r_1>0$ such that for all $k\geq k_1$, for all $x\in \{x\in \R^n: \psi_{\gk}(x)=0\}$, and for all $z\in B_{r_1}(x)$, we have 
 \begin{equation*}\label{psiuniform}
  \|\nabla\psi_{\gk}(z)\|> 2\eta.
\end{equation*} 
In particular,  for $k\geq k_1$ we have \begin{equation}\label{A2.2forpsigk} [\psi_{\gk}(x)=0]\im \|\nabla\psi_{\gk}(x)\|>2\eta.\end{equation}
\item There exists $k_2\geq k_1$ and $\e_{o}>0$ such that for all $k\geq k_2 $ we have  \begin{equation*}\label{epsilonpsigk} \big[x\in C^{\gk}\;\hbox{and}\; \|\nabla\psi_{\gk}(x)\|\leq \eta\big]\implies \psi_{\gk}(x)<-\e_o.\end{equation*}
\end{enumerate}
\end{proposition}

\begin{remark} \label{xi.bounded.onCk}  One can easily see from the definition of $\psi_{\gk}$ in \eqref{psigkdef} that when $r>1$, the second inequality of \eqref{pisgkineq} is strict, that is,  $\psi_i(x)< \psi_{\gk}(x)$ for all $x\in\R^n$ and for $i=1,\dots,r$. On the other hand, from \eqref{sigmadef}$(a)$, \eqref{Cgkkdef} and \eqref{pisgkineq}, we deduce that for all $k\in \N$ and for $i=1,\dots,r$, we have 
$$\gk e^{\gk\psi_i(x)}\leq \gk e^{\gk\psi(x)}\leq \gk e^{\gk\psi_{\gk}(x)}\leq  \frac{2M}{\eta},\,\;\forall x\in C^{\gk}(k).$$

\end{remark}

Employing  the preceding proposition and \cite[Proposition 3.1]{verachadi}, we show the following properties for the sequence of sets $(C^{\gk})_k$. One novelty of this Proposition is provided in the second part  of its item $(c)$, namely, that for  $k$ large enough, the sets  $C^{\gk}$ are {\it uniformly} prox-regular with a uniform prox-regularity constant being   $\frac{\eta}{2M_{\psi}} $. Note that 
in  \cite[Proposition 3.1]{verachadi} it is only established that  $C^{\gk}$ is $\frac{\eta}{M_{\psi_{\gk}}} $-prox-regular, where $2M_{\psi_{\gk}}$ depends on $k$, for being the Lipschitz constant of $\nabla \psi_{\gk}$ over the compact set $\conv C$.   Thus,  establishing the {\it uniform} prox-regularity is  not straightforward, since    the generalized Hessian of $\psi_{\gk}$ is {\it unbounded} and so is the sequence $(M_{\psi_{\gk}})_k$.

\begin{proposition} \label{propcgk} For all $k$, the set  $C^{\gk}\subset C$ is compact with $C^{\gk}\subset \inte C$ for $r>1$. Moreover, there exists $k_3\geq k_2$ such that for $k\geq k_3$, we have$\sp:$
\begin{enumerate}[$(a)$]
\item $\bdry C^{\gk}=\{x\in\R^n : \psi_{\gk}(x)=0\}=\left\{x\in\R^n : \sum_{i=1}^r e^{\gk\psi_i(x)}=1\right\}\not=\emptyset.$
\item $\inte C^{\gk}=\{x\in\R^n : \psi_{\gk}(x)<0\}\not=\emptyset.$
\item $C^{\gk}$ is epi-Lipschitzian, $C^{\gk}=\clo(\inte C^{\gk})$, and  $C^{\gk}$ is $\frac{\eta}{2M_{\psi}} $-prox-regular.
\item For all $x\in\bdry C^{\gk}$ we have \begin{eqnarray*} N_{C^{\gk}}(x)=N^P_{C^{\gk}}(x)=N^L_{C^{\gk}}(x)&=&\bigg\{\l\nabla\psi_{\gk}(x) : \l\geq 0\bigg\}\\&=&\bigg\{\l\sum\limits_{i=1}^r e^{\gk\psi_i(x)}\nabla\psi_i(x) : \l\geq 0\bigg\}.\end{eqnarray*}
\end{enumerate}
Furthermore, the sequence $(C^{\gk})_k$ is a nondecreasing sequence whose  Painlev\'e-Kuratowski limit is $C$ and satisfies \begin{equation}\label{union.Cgk}\inte C=\bigcup_{k\in \N} \inte C^{\gk}=\bigcup_{k\in \N} C^{\gk}. \end{equation}
\end{proposition}
{\it Proof}  From \eqref{pisgkineq} and the definition of $C^{\gk}$ in \eqref{Cgkdef}, we conclude that $C^{\gk}\subset C$ for all $k$. This gives that the closed set $C^{\gk}$ is bounded for all $k$, and hence  $C^{\gk}$ is compact for all $k$. On the other hand, for $r>1$, if  $x\in C^{\gk}$ then by \eqref{Cgkdef}, we have $\psi_{\gk}(x)\le 0,$  which yields that $\psi_i(x)<0$ for  $i=1,\dots,r$, by Remark \ref{xi.bounded.onCk}. Hence, $$x\in\bigcap_{i=1}^{r} \inte C_i=\inte C.$$ 
Therefore, for $r>1$, we have $C^{\gk}\subset \inte C$ for all $k$.

As $\inte C\not=\emptyset,$ then,  for any $x\in\inte C$  we have $x\in C$ and  $\psi(x)<0.$ Using \eqref{pisgkineq}, there exists $k_3\geq k_2$, such that  for $k\geq k_3$, we have that $$\psi_{\gk}(x)\leq\psi(x)+\frac{\ln(r)}{\gk}<0.$$
This gives that $x\in C^{\gk}$, and hence $C^{\gk}\not=\emptyset$ for  $k\geq k_3$. Note that both arguments made  above yield that  
$$\inte C\subset \bigcup_{k\in \N} \inte C^{\gk}  \subset \bigcup_{k\in \N} \inte C^{\gk}  \subset \inte C.$$ This gives that  \eqref{union.Cgk} holds true.
Hence, since by Proposition \ref{proppsigk},  $\psi_{\gk}$ is $C^{1,1}$ and satisfies \eqref{A2.2forpsigk}, we deduce that all the properties satisfied by the set $C$ in \cite[Proposition 3.1]{verachadi} are also satisfied by $C^{\gk}$ for all $k\geq k_2$. Therefore, the assertions $(a)$-$(d)$ of Proposition \ref{propcgk} are valid except the {\it uniform} constant for the  prox-regularity of $C^{\gk}$.  For, let $x\in\bdry C^{\gk}$ and let $0\not=\zeta\in N^P_{C^{\gk}}(x)$. Then we have, for some $\l>0$, that \begin{equation}\label{verynice}\zeta=\l\sum\limits_{i=1}^r e^{\gk\psi_i(x)}\nabla\psi_i(x)\;\;\hbox{and}\;\;\sum\limits_{i=1}^r e^{\gk\psi_i(x)}=1.\end{equation}
For $y$ fixed in $C^{\gk}$ and for $t\in (0,1]$, we have $x+t(y-x)\in\conv C^{\gk}\subset \conv C$. Since  $2M_\psi$ is a common Lipschitz constant of the finite family $\{\nabla \psi_{i}\}_{i=1}^{r}$, it follows that,  for $i=1,\dots,r$, 
\begin{eqnarray*} \<\nabla\psi_i(x+t(y-x))\bp-\bp\nabla\psi_i(x),y-x\>&=& \frac{1}{t} \<\nabla\psi_i(x+t(y-x))\bp-\bp\nabla\psi_i(x),t(x-y)\> \\* &\geq&-2tM_{\psi}\|x-y\|^2.\end{eqnarray*} Hence, using the mean value theorem, we have for $i=1,\dots,r$ the existence of $t^i_y\in(0,1)$ such that  \begin{eqnarray*}0\geq \psi_i(y)-\psi_i(x)&=&\<\nabla\psi_i(x+t^i_y(y-x)), y-x\> \\&\geq& -2t^i_yM_{\psi}\|x-y\|^2 + \<\nabla\psi_i(x),y-x\>.\end{eqnarray*}
Whence, for $i=1,\dots,r,$ $$ \<\nabla\psi_i(x),y-x\>\leq 2t^i_yM_{\psi}\|y-x\|^2\leq 2M_{\psi}\|y-x\|^2. $$
 Using \eqref{verynice}, this gives  that $$\<\zeta,y-x\>\leq 2\l M_{\psi}\|y-x\|^2.$$
Therefore, from \eqref{gradpsigk}, \eqref{A2.2forpsigk} and \eqref{verynice}, we conclude that $$\left\<\frac{\zeta}{\|\zeta\|},y-x\right\>\leq \frac{2\l M_{\psi}}{\|\zeta\|}\|y-x\|^2=\frac{2\l M_{\psi}}{\l\|\nabla\psi_{\gk}(x)\|}\|y-x\|^2\leq \frac{M_{\psi}}{\eta}\|y-x\|^2.$$
This terminates the proof of the $\frac{\eta}{2M_{\psi}}$-prox-regularity of $C^{\gk}$. 

We proceed to prove the \enquote{Furthermore} part of Proposition \ref{propcgk}. Since $(\psi_{\gk})_k$ is monotonically nonincreasing  in $k$, we deduce that the sequence $(C^{\gk})_k$ is a nondecreasing sequence. Hence, it is easy to show (see e.g., \cite[ Exercice 4.3]{rockwet}) that the Painlev\'e-Kuratowski limit of the sequence $(C^{\gk})_k$ satisfies 
\begin{equation}\label{limitCgk} \lim_{k\f\infty} C^{\gk}=\clo\Bigg(\bigcup_{k\in \N} C^{\gk}\Bigg). \end{equation}
Now, upon taking the closure of $\inte C$ in the already established \eqref{union.Cgk} and using from Proposition \ref{prop1} that $C=\clo({\inte C})$, equation \eqref{limitCgk} yields that the Painlev\'e-Kuratowski limit of the sequence $(C^{\gk})_k$ is $C$. \qed

We proceed to present the properties of the sequence of sets $(C^{\gk}(k))_k$. For $c\in \bdry C$, we denote by $d_c$ the nonzero vector $d_c:=\sum_{j\in\I^0_c}v_j(c) $, where for $j=1,\dots,r$, $v_j(c)$ is the unique projection of $-\nabla\psi_j(c)$ to the Clarke tangent cone $T_C(c)$. For more information about the vector $d_c$, see \cite[Lemma 6.1]{verachadijca}. Note that when $r=1$, the vector $d_c$ coincides with  $-\nabla\psi_j(c)\in T_C(c)$. 

\begin{proposition}[{\cite[Proposition 4.3 \& Remark 4.4]{verachadijca}}]\label{propcgk(k)} 
The following assertions hold\sp$:$
\begin{enumerate}[$(i)$, leftmargin=\widthof{(iii)}+\labelsep]
\item  For all $k$, the set $C^{\gk}(k)\subset \inte C^{\gk}\subset \inte C$ and is compact. Moreover, there exists $k_4\geq k_3$ such that for $k\geq k_4$, we have$\sp:$
\begin{enumerate}[$(a)$]
\item $\bdry C^{\gk}(k)=\{x\in\R^n : \psi_{\gk}(x)=-\a_k\}\\\hspace*{1.85cm}=\left\{x\in\R^n : \sum_{i=1}^r e^{\gk\psi_i(x)}= \frac{2{M}}{\eta\gk}\right\}\not=\emptyset.$
\item $\inte C^{\gk}(k)=\{x\in\R^n : \psi_{\gk}(x)<-\a_k\}\not=\emptyset.$
\item $[\psi_{\gk}(x)=-\a_k]\im \|\nabla\psi_{\gk}(x)\|>\eta.$
\item $C^{\gk}(k)$ is $\frac{\eta}{4M_{\psi}}$-prox-regular,\textnormal{\footnote{In \cite[Proposition 4.3]{verachadijca}, the prox-regularity constant of the set $C^{\gk}(k)$ was obtained to be $\frac{\eta}{2M_{\psi_{\gk}}}$, where $2M_{\psi_{\gk}}$ is the Lipschitz constant of $\nabla \psi_{\gk}$ over the compact set $\conv C$. Using arguments similar to those used in the proof of Proposition \ref{propcgk}, one can prove that $M_{\psi_{\gk}}$ can be replaced by $2M_{\psi}$.}}  and epi-Lipschitz with $$C^{\gk}(k)=\clo(\inte C^{\gk}(k)).$$
\item For all $x\in\bdry C^{\gk}(k)$ we have \begin{eqnarray*}N_{C^{\gk}(k)}(x)\bbp= \bbp N^P_{C^{\gk}(k)}(x)&=&N^L_{C^{\gk}(k)}(x))\\&=&\bigg\{\l\nabla\psi_{\gk}(x) : \l\geq 0\bigg\}\\&=&\bigg\{\l\sum\limits_{i=1}^r e^{\gk\psi_i(x)}\nabla\psi_i(x) : \l\geq 0\bigg\}.\end{eqnarray*}
\end{enumerate}
\item  The sequence $(C^{\gk}(k))_k$ is a nondecreasing sequence whose  Painlev\'e-Kuratowski limit is $C$ and satisfies
\begin{equation}\label{union.Cgk(k)}\inte C=\bigcup_{k\in \N} \inte C^{\gk}(k)=\bigcup_{k\in \N} C^{\gk}(k). \end{equation}
\item  For $c\in \bdry C$, there exist ${k}_c\geq k_4$ and ${r}_{\c}>0$ such that  \begin{equation*}\label{lastideahope}\left(\left[C\cap \bar{B}_{{r}_{\c}}({c})\right]+\sigma_k\frac{d_c}{\|d_c\|}\right)\subset \inte C^{\gk}(k)\subset \inte C^{\gk},\;\;\forall k\geq {k}_c.\end{equation*}
In particular,  for $k\geq {k}_c$ we have \begin{equation*}\label{boundarypoint(k)}
 \left({c}+\sigma_k\frac{d_c}{\|d_c\|}\right)\in \inte C^{\gk}(k)\subset\inte C^{\gk}. \end{equation*}
\end{enumerate}
\end{proposition}

For $x_0$ the initial point of the problem $(P)$, we define the sequence $\left(x_0^k\right)_k$ by
\begin{equation} \label{initial} x_0^k:= 
     \begin{cases}x_0,\;\forall k\in\N, &\;\hbox{if}\;x_0\in\inte {C},\vspace{0.2cm}\\ \displaystyle {x_0}+\sigma_k \frac{d_{x_0}}{\|d_{x_0}\|},\;\forall k\in\N,&\;\hbox{if}\; x_0\in\bdry C.\end{cases} \end{equation}
Since $\sigma_k\f 0$, the following lemma follows from \eqref{union.Cgk(k)} and Proposition \ref{propcgk(k)}$(iii)$.

\begin{lemma}\label{lem1} The sequence $x_0^k$ converges to $x_0$, and there exists $k_5\geq k_{x_0}$ such that $x_0^k\in C^{\gk}(k)$ for all $k\geq k_5$.
\end{lemma}

\begin{remark} From (H3) we can deduce that for $i=1,\dots,r$, the set $C_i$ satisfies the same assumptions satisfied by $C$ in the papers \cite{VCpaper,verachadinum,verachadi}. Hence, all the properties established in those papers for $C$ and $C(k)$ are valid here for $C_i$ and $C_i(k)$, respectively, where $i=1,\dots,r$. On the other hand, from \eqref{pisgkineq} and Remark \ref{xi.bounded.onCk}, we have \begin{equation} \label{cgk(k)toci(k)} C^{\gk}(k)\subset \bigcap_{i=1}^{r} C_i(k),\end{equation}
and when $r>1$, $C_i(k)$ can be replaced in \eqref{cgk(k)toci(k)} by $\inte  C_i(k)$.
\end{remark}

We terminate this section with the following proposition in which we provide properties of the projection maps from $\bdry C$  to $\bdry C^{\gk}(k)$ and from $\bdry C^{\gk}$  to $\bdry C^{\gk}(k)$.

\begin{proposition} There exists  $k_6\geq k_5$ such that for $k\geq k_6$, the projection map $\pi_k\colon \bdry C\f \bdry C^{\gk}(k)$ is onto and 2-Lipschitz, and the projection map $\hat{\pi}_k\colon \bdry C^{\gk}\f \bdry C^{\gk}(k)$ is bi-Lipschitz.
\end{proposition}
{\it Proof} We begin with the projection map $\pi_k$. Since $C$ is compact, and, by Proposition \ref{propcgk(k)}, $C$ is the Painlev\'e-Kuratowski limit of $C^{\gk}(k)$, we deduce that $d_H(C,C^{\gk}(k))\f 0$ as $k\f 0$. This gives that, for $k$ sufficiently large,  $d(c,{\bdry C^{\gk}}(k))<\frac{\eta}{4M_\psi}$ for all $c\in \bdry C$. In addition, by Proposition \ref{propcgk(k)}$(d)$, we have that $C^{\gk}(k)$ is $\frac{\eta}{4M_\psi}$-prox-regular for $k$  large enough. We conclude that, for $k$ sufficiently large, $\pi_k$ is a single valued function. Now, by taking $k$ large enough so that $d(c,{\bdry C^{\gk}}(k))<\frac{\eta}{8M_\psi}$ for all $c\in \bdry C$, we deduce from \cite[Theorem 4.8]{prox} that $\pi_k$ is 2-Lipschitz. We claim that ${k}$ can be taken large enough so that $\pi_k$ is an onto function. Indeed, due to the $\frac{\eta}{4M_\psi}$-prox-regularity of $C^{\gk}(k)$,  it is sufficient to prove that for $k$ large enough, we have for each $c\in\bdry C^{\gk}(k)$, $$(\bdry C)\cap \left[c,c+\frac{\eta}{4M_\psi}\frac{\nabla\psi_{\gk}({c})}{\|\nabla\psi_{\gk}({c})\|}\right)\not=\emptyset.$$ If not, then there exist an increasing sequence $k_n\in\N$ and a sequence $c_n\in \bdry C^{\gkn}(k_n)$ such that $$(\bdry C)\cap \left[c_n,c_n+\frac{\eta}{4M_\psi}\frac{\nabla\psi_{\gkn}({c_n})}{\|\nabla\psi_{\gkn}({c_n})\|}\right)=\emptyset.$$ This gives that $$\psi_{\gkn}(c_n)=-\a_{k_n}\;\,\hbox{and}\;\,\psi\left(c_n+\frac{\eta}{8M_\psi}\frac{\nabla\psi_{\gkn}({c_n})}{\|\nabla\psi_{\gkn}({c_n})\|}\right)<0,\;\;\forall n.$$  
Using the compactness of $C$, the convergence of $\a_{k_n}$ to $0$, and the inequalities of \eqref{pisgkineq},  it follows that there exist  a subsequence of $k_n$, we do not relabel, a $\bar{c}\in C$, and a unit vector $\bar{v}$, such that 
$$c_n\f\bar{c}\in\bdry C,\;\, \frac{\nabla\psi_{\gkn}({c_n})}{\|\nabla\psi_{\gkn}({c_n})\|}\f \bar{v},\;\,\hbox{and}\;\,\psi\left(\bar{c}+\frac{\eta}{8M_{\psi}}\bar{v}\right)\leq 0.$$ 
Hence, \begin{equation}\label{pionto} \bar{c}+\frac{\eta}{8M_{\psi}}\bar{v}\in C.\end{equation} 
Since  each point in $\bdry C$ is the limit of a sequence of points in $C^{\gkn}(k_n)$,  then, the  $\frac{\eta}{4M_{\psi}}$-prox-regularity of $C^{\gkn}(k_n)$ implies that  $\bar{v}\in N_C^P(\bar{c})$.   Thus, the $\frac{\eta}{M_{\psi}}$-prox-regularity of $C$ yields that $$\bar{c}+\frac{\eta}{8M_{\psi}}\bar{v}\in B_{\frac{\eta}{M_{\psi}}}\left(\bar{c}+\frac{\eta}{M_{\psi}}\bar{v}\right)\subset C^c,$$ which contradicts \eqref{pionto}.

We proceed to prove that $\hat{\pi}_k$ is bi-Lipschitz. By Proposition \ref{propcgk}, the function $\psi_{\gk}$ satisfies the same assumptions satisfied by the function $\psi$ in \cite{VCpaper}. In addition,  the two sets $C^{\gk}$ and $C^{\gk}(k)$ are  defined in terms of  $\psi_{\gk}$  in the same way $C$ and $C(k)$ of \cite{VCpaper} were defined in terms of $\psi$. Hence, from \cite[Theorem 3.1$(vi)$]{VCpaper}, we can deduce that, for $k$ large enough, $\hat{\pi}_k$ is bi-Lipschitz. \qed

\begin{remark} \label{piknotbi} Unlike the case $r=1$ studied in \cite[Theorem 3.1$(vi)$]{VCpaper} and   \cite[Lemma 2]{verachadinum}, the projection  $\hat{\pi}_k$ here for $r>1$ cannot be shown to be  {\it uniformly} bi-Lipschitz. The reason is that, when $r>1$, the generalized Hessian of the function $\psi_{\gk}$ is {\it not} uniformly bounded.  More issues surface for the projection ${\pi}_k$, since the function $\psi$ is only Lipschitz.  These facts render the techniques used in the proof of \cite[Theorem 3.1$(vi)$]{VCpaper} not applicable for either projections.

\end{remark}

\section{Key Results} \label{key} Parallel to \cite[Section III]{verachadinum}, we provide in this section three theoretical results, that are the keystone of our numerical algorithm constructed for $(P)$. We note that having  $\psi$  {\it only} Lipschitz and the generalized Hessian of $\psi_{\gk}$ {\it not} bounded in the general case $r>1$,  make the proofs of these results more challenging than their counterparts for the case $r=1$, and hence, new ideas and techniques are required.

For $x^k_0$ given in \eqref{initial}, we denote by $(D_{\gk})$  the approximation dynamic defined by  
\begin{equation}\hspace{-0.2cm}\label{Dgk1} ({D}_{\gk}) \begin{cases} \dot{x}(t)=f(x(t),u(t))-\sum\limits_{i=1}^{r} \gk e^{\gk\psi_i(x(t))} \nabla\psi_i(x(t)),\;\textnormal{a.e.}\; t\in[0,T],\\ x(0)=x^k_{0}. \end{cases}
 \end{equation}
One can easily verify  that using \eqref{gradpsigk},  the system  $(D_{\gk})$  can be rewritten in terms of $\psi_{\gk}$ as  follows:
 \begin{equation}\hspace{-0.25cm} \label{Dgk2} ({D}_{\gk}) \begin{cases} \dot{x}(t)=f(x(t),u(t))-\gk e^{\gk\psi_{\gk}(x(t))} \nabla\psi_{\gk}(x(t)),\;\textnormal{a.e.}\; t\in[0,T],\\ x(0)=x_0^{k}. \end{cases} \end{equation}
For a solution $x_{\gk}$  of $(D_{\gk})$ corresponding to a control $u\in \U$, we denote by $\xi_{\gk}(\cdot)$ the sequence of non-negative continuous functions on $[0,T]$ defined by
 \begin{equation*} \label{defxi} \xi_{\gk}(\cdot):= \gk e^{\gk\psi_{\gk}(x_{\gk}(\cdot))} \overset{\eqref{psigkdef}}{=}\sum_{i=1}^{r}\xi_{\gk}^i(\cdot),\end{equation*} where $\xi_{\gk}^i(\cdot):=\gk e^{\gk\psi_i(x_{\gk}(\cdot))}$ for $i=1,\dots,r$.
 
From \cite[Theorem 4.13]{verachadijca}, we can deduce the following proposition.
\begin{proposition}  \label{propo1} There exists $k_7\geq k_6$ such that for all $k\geq k_7$ and for all $u\in U$, the solution $x_{\gk}$ of $(D_{\gk})$ corresponding to $u$ satisfies\sp$:$
\begin{enumerate}[$(i)$, leftmargin=\widthof{(iii)}+\labelsep]
\item $x_{\gk}(t)\in C^{\gk}(k)\subset \inte C^{\gk}\subset \inte C$ for all $t\in[0,T]$.
\item $0< \xi_{\gk}^i(t) \leq \xi_{\gk}(t)\le \frac{2{M}}{\eta}$ for all $t\in [0,T]$ and for $i=1,\dots,r.$
\item $ \|\dot{x}_{\gk}(t)\|\leq {M}+\frac{2{M}\bar{M}_\psi}{\eta}$ \,for \textnormal{a.e.} $t\in [0,T]$.
\end{enumerate}
\end{proposition}

In the following proposition, we prove that the $L^{\infty}$-distance between the solution $x_{\gk}$ of $(D_{\gk})$ and the solution $x$ of $(D)$  is controlled by $\sigma_k$, when  the same control $u\in \U$ is used in both  dynamics. Note that this result {\it cannot} be  deduced from  \cite[Proposition 2]{verachadinum}, where the  function $\psi$ is replaced by our function $\psi_{\gk}$. The reason  behind this is that unlike the case for $\psi$ in  \cite[Proposition 2]{verachadinum}, the Lipschitz constant for  $\nabla\psi_{\gk}$ is {\it not uniformly} bounded for $k$ large.

\begin{proposition}  \label{propo2} 
There exists $k_8\geq k_7$ such that for all $k\geq k_8$ and for all $u\in \mathscr{U}$, the solution $x$ of system $(D)$ and the solution $x_{\gk}$ of system $(D_{\gk})$, both corresponding to the same control $u\in\U$, satisfy 
\begin{equation*}\label{estimate} \|x_{\gk}(t)-x(t)\|^2\leq e^{\tilde{M}T} \sigma_k^2 + \frac{8\eta M (e^{\tilde{M}T}-1)}{\tilde{M}M_{\psi}}\sigma_k,\;\,\forall t\in [0,T]\;\textnormal{a.e.}, \end{equation*}
where $\tilde{M}:=\frac{5MM_\psi}{ \eta}+2{M}$.
\end{proposition}
{\it Proof} As the Lipschitz constant of $\nabla\psi_{\gk}$ is {\it not} uniformly bounded in $k$, a modification of the proof of \cite[Proposition 2]{verachadinum} is required here. For this, we shall use the version \eqref{Dgk1}  of $(D_{\gk})$   instead of that in \eqref{Dgk2}.   Now,  given that for $i=1,\dots,r$, $\psi_i$ is $\CO^{1,1}$,  the second-order generalized Taylor expansion,  \cite[Theorem 2.3]{hiriart}, implies that for  $x$ and $y$ in $C$, there exist $z_i\in \conv (C)$ and $M_{z_i}\in \partial^2\psi_i(z_i)$ such that for $i=1,\dots,r$, 
\begin{equation*} \label{hiriart}\psi_i(y)=\psi_i(x)+\<\nabla \psi_i(x),y-x\>+\frac{1}{2}\<M_{z_i}(y-x),y-x\>.\end{equation*}
Hence, using that for $i=1,\dots,r$, $\|M_{z_i}\|\le 2 M_\psi$ and  $\psi_i(y)\le 0$, we obtain that for all $x$ and $y$ in $C$,
\begin{equation} \label{grad}
-\<\nabla \psi_i(x),x-y\>\le -\psi_i(x)+M_\psi\|x-y\|^2\;\;\hbox{for}\;i=1,\dots,r,
\end{equation}
Employing (H1), (\ref{grad}), the inclusion of $(D)$, the prox-regularity of $C$,  the version \eqref{Dgk1} of $(D_{\gk})$, Lemma \ref{lem1}, and the uniform boundedness of $(\xi^i_{\gk})_{i=1}^r$ in Proposition \ref{prop1}, we obtain that, for $k\geq k_7$ and for a.e. $t\in[0,T]$,  
\begin{eqnarray} \nonumber & &\hspace{-0.3cm}\frac{d}{dt}\|x_{\gk}(t)-x(t)\|^2
= 2\<x_{\gk}(t)-x(t),\dot{x}_{\gk}(t)-\dot{x}(t)\>\\
&=& \nonumber2\<x_{\gk}(t)-x(t),f(x(t),u(t))-\dot{x}(t)\>\\
&+& \nonumber 2\<x_{\gk}(t)-x(t),f(x_{\gk}(t),u(t))-f(x(t),u(t))\>\\
& \nonumber-&2\sum_{i=1}^r\xi^i_{\gk}(t)\<x_{\gk}(t)-x(t),\nabla\psi_i(x_{\gk}(t))\>\\
&\leq & \nonumber \frac{MM_\psi}{\eta}\left\|x_{\gk}(t)-x(t)\right\|^2 + 2M\left\|x_{\gk}(t)-x(t)\right\|^2 \\&-&\nonumber2\sum_{i=1}^r\xi^i_{\gk}(t)\psi_i(x_{\gk}(t))+\nonumber2 M_\psi \left\|x_{\gk}(t)-x(t)\right\|^2\sum_{i=1}^r \xi^i_{\gk}(t)\\ &\leq &\nonumber \frac{MM_\psi}{\eta}\left\|x_{\gk}(t)-x(t)\right\|^2 + 2M\left\|x_{\gk}(t)-x(t)\right\|^2\\&-& 2\sum_{i=1}^r \xi^i_{\gk}(t)\psi_i(x_{\gk}(t))+
\nonumber\frac{4MM_\psi}{\eta}\left\|x_{\gk}(t)-x(t)\right\|^2 \\&=&\tilde{M} \label{distanceest} \left\|x_{\gk}(t)-x(t)\right\|^2-2\sum_{i=1}^r \xi^i_{\gk}(t)\psi_i(x_{\gk}(t)),\end{eqnarray} 
where $\tilde{M}=\frac{5MM_\psi}{ \eta} +2 M$. Choose $k_8\geq k_7$ sufficiently large, so that $\gk\geq \frac{2Me}{\eta}$ for all $k\ge k_8.$ This means that $a:=\frac{2M}{\eta\gk}=e^{-\gk \a_k}\le e^{-1}$ for all $k\ge k_8$.
Then, using the facts that $xa^x\leq a $ for $a\leq e^{-1}$ and $x\geq 1$,   and that  $\psi_i(x_{\gk}(t))\leq\psi(x_{\gk}(t))\leq \psi_{\gk}(x_{\gk}(t))  \leq -\a_k$ for $i=1,\dots,r$ and for all $k\geq {k}_8$, we deduce with the help of \eqref{sigmadef}$(a)$ that   for $i=1,\dots,r$,
\begin{eqnarray*} -\xi^i_{\gk}(t)\psi_i(x_{\gk}(t))&=& (\a_k\gk) \left(\frac{2M}{\eta\gk}\right)^{\left(\frac{-\psi_i(x_{\gk}(t))}{\a_k}\right)}\left(\frac{-\psi_i(x_{\gk}(t))}{\a_k}\right)\\[4pt]&\leq& (\a_k\gk)\, \frac{2M}{\eta\gk} = \frac{2M\a_k}{\eta}, \quad \forall\sp t\in[0,T],\;\; \forall \sp k\ge {k}_1.
\end{eqnarray*}
Hence by (\ref{distanceest}) we conclude that for $k\ge k_8$  and for a.e. $t\in[0,T]$ $$ \frac{d}{dt}\|x_{\gk}(t)-x(t)\|^2\leq \tilde{M}  \left\|x_{\gk}(t)-x(t)\right\|^2+\frac{4rM\a_k}{\eta}.$$
Now using Gronwall's lemma \cite[Proposition 4.1.4]{clsw}, the definition of $\sigma_k$ in \eqref{sigmadef}, and (\ref{initial}), we get that for $k\ge {k}_8$ and for all\;$t\in[0,T]$, 
\begin{eqnarray*}\|x_{\gk}(t)-x(t)\|^2&\leq & e^{\tilde{M}T} \|x_0^k-x_0\|^2+ \frac{4rM(e^{\tilde{M}T}-1)}{\eta\tilde{M}}\a_k\\& \leq& e^{\tilde{M}T} \sigma_k^2+ \frac{4r{M}(e^{\tilde{M}T}-1)}{\eta\tilde{M}}\a_k \\ &\leq& e^{\tilde{M}T} \sigma_k^2 + \frac{8\eta M (e^{\tilde{M}T}-1)}{\tilde{M}M_{\psi}}\sigma_k.\end{eqnarray*}
This terminates the proof of the proposition. \qed

The following lemma is a generalization to the case when $r>1$ of  \cite[Lemma 2]{verachadinum}, established for  $r=1$, and hence, our $C^{\gk}(k)$ here reduces to  $C(k)$ there.  Note that  the proof of  \cite[Lemma 2]{verachadinum} is based on the {\it uniform bi-Lipschitz} continuity of the projection map from $\bdry C$ to $\bdry C(k)$;  a property {\it not} met here for the projection map from $\bdry C$ to $\bdry C^{\gk}(k)$, see Remark \ref{piknotbi}. Therefore,  new ideas are needed here in order to prove this lemma. One may believe that this result could be established by replacing $C$  with $C^{\gk}$  in that projection, that is, to use  the projection map $\hat{\pi}_k$ from $\bdry C^{\gk}$ to $\bdry C^{\gk}(k)$. But, as mentioned in Remark \ref{piknotbi}, this projection map is {\it not} necessarily uniformly bi-Lipschitz. Moreover,  the quasiconvexity  required in the proof of \cite[Lemma 2]{verachadinum} for $\bdry C$, is {\it not} guaranteed to be {\it uniform} here for $\bdry C^{\gk}.$

\begin{lemma} \label{lem2} There exist $k_9\geq k_8$ and $\beta>0$ such that for all $k\geq k_9$, $$\<\nabla e^{\gk\psi_{\gk}(y)}- \nabla e^{\gk \psi_{\gk}(x)},y-x\>\ge -\frac{4r\beta M M_{\psi}}{\eta}\|y -x\|^2,\;\,\forall x,\sp y\in C^{\gk}(k).$$
\end{lemma}
{\it Proof} From \cite[Lemma 2]{verachadinum} applied to each $C_i$ for $i=1,\dots,r$, we get the existence of $k^i_9\geq k_8$ and $\beta_i>0$ such that for all $k\geq k^i_9$

$$\<\nabla e^{\gk\psi_i(y)}- \nabla e^{\gk \psi_i(x)},y-x\>\ge -\frac{4\beta_i M M_{\psi}}{\eta}\|y -x\|^2,\;\,\forall x,\sp y\in C_i(k).$$
Using that $e^{\gk\psi_{\gk}(\cdot)}=\sum_{i=1}^re^{\gk\psi_i(\cdot)}$ and the inclusion \eqref{cgk(k)toci(k)}, we conclude that for $k_9:=\max\{k_9^i : i=1,\dots,r\}$ and $\beta:=\max\{\beta_i : i=1,\dots,r\}$, we have for all $k\geq k_9$, $$\<\nabla e^{\gk\psi_{\gk}(y)}- \nabla e^{\gk \psi_{\gk}(x)},y-x\>\ge -\frac{4r\beta M M_{\psi}}{\eta}\|y -x\|^2,\;\,\forall x,\sp y\in C^{\gk}(k).$$
The proof of the lemma is terminated.\qed

\begin{remark}\label{convex} When for all $i=1,\dots,r$, the set $C_i$ is convex (that is, $\psi_i$ is convex), then the result of Lemma \ref{lem2} becomes straightforward since in that case the function $\psi_{\gk}$ is convex, see \cite[Lemma 3]{LiFang}, which gives that $e^{\gk\psi_{\gk}(\cdot)}$ is convex, and hence, $$\<\nabla e^{\gk\psi_{\gk}(y)}- \nabla e^{\gk \psi_{\gk}(x)},y-x\>\geq 0,\;\,\forall x,\sp y\in\R^n.$$
Note that for $n>1$, the convexity of $C_i$ implies that $\bdry C_i$ is connected.
\end{remark}

Let $N$ be a positive integer and set $h:=\frac{T}{N}$. For any vector $\big(u_1^N,\dots, u_{N}^N)$ in $U^N:=\underbrace{U\times\cdots\times U}_{N}$, we associate the piecewise constant control function $$u^N(t)\equiv u_j^N\;\,\hbox{for}\;\,t\in [(j-1)h,jh)\;\,\hbox{and}\; j=1,\dots,N.$$ We denote by $\U^N$ the set of such controls, and by $x^N_{\gk}$ the solution of $(D_{\gk})$  corresponding to a control $u^N\in\U^N$.

As in \cite[Section III]{verachadinum}, Lemma \ref{lem2} leads to the following proposition, whose proof follows arguments similar to those used in the proof of \cite[Proposition 3]{verachadinum}.

\begin{proposition} \label{propo3} Let $u\in \mathscr{U}$ and, for $k\ge {k}_9$, let $x_{\gk}$ be the solution of $(D_{\gk})$ corresponding to $u$. Then there exists  $u^N\in\U^N$ such that $x^N_{\gk}$, the solution of $(D_{\gk})$ corresponding 
to  $u^N,$ satisfies for $\delta_j:=\|x_{\gk}(jh)-x^N_{\gk}(jh)\|^2$ the inequality  
\begin{equation*} \label{xN-xgkN}
\delta_j\leq \frac{2\hat{M}^2 e^{6\bar{M}T}}{3\bar{M}}  (1+6 \bar{M} h)(1+\bar{M}h)h,\,\;\;\hbox{for}\;\,j=1,\dots,N,\end{equation*}
where $\bar{M}:=M+\frac{4r\beta MM_{\psi}}{\eta}$ and $\hat{M}:={M}+\frac{2{M}\bar{M}_\psi}{\eta}.$ \end{proposition}

\section{Numerical Algorithm} \label{algo} Based on the three key results of the previous section, namely Propositions \ref{propo1}, \ref{propo2} and \ref{propo3}, we prove in this section the main result of this paper which inspires  the construction of  our numerical algorithm that solves $(P)$. We denote by $(P_{\gk})$ the problem $(P)$ in which the dynamic $(D)$ is replaced by $(D_{\gk})$, that is, $$ \begin{array}{l} (P_{\gk})\colon\; \hbox{Minimize}\;g(x(T))\\ \hspace{1.17cm} \hbox{over}\;(x,u)\;\hbox{such that}\;u\in \mathscr{U}, x\in AC([0,T];\R^n), \;\hbox{and} \\[2pt]  \hspace{1.17cm} ({D}_{\gk}) \begin{cases} \dot{x}(t)=f(x(t),u(t))-\gk e^{\gk\psi_{\gk}(x(t))} \nabla\psi_{\gk}(x(t)),\,\;\textnormal{a.e.}\; t\in[0,T],\\ x(0)=x_0^{k}.  \end{cases}
    \end{array}$$

We fix $\e>0$. From Propositions \ref{propo2} and \ref{propo3}, we have the following:
\begin{enumerate}[$\bullet$]
\item Since $\sigma_k\f 0$ in Proposition \ref{propo2}, then, there exists $k_{\e}\ge k_9$ such that  for a given solution $(\x,\u)$ of $(P)$ and  for $\x_{\gk}$ being the solution of $(D_{\gk})$ corresponding to $\u$, we have  \begin{equation*} \label{num1} \|\bar{x}-\x_{\gk}\|_{\infty}\leq \frac{\e}{3L_g},\;\;\;\forall k\ge k_{\e}.\end{equation*} 
\item From Proposition \ref{propo3}, we deduce the existence of a positive integer $N_{\e}$ such that for $N\ge N_{\e}$ and $h:=\frac{T}{N}$, there is $\u^N\in\U^N$ for which the solution $\x^N_{\gk}$ of $(D_{\gk})$ corresponding to $\u^N$ satisfies  \begin{equation*}  \label{num2} \|\x_{\gk}-\x_{\gk}^N\|_{\infty}\leq \frac{\e}{3L_g},\;\;\;\forall k\ge k_{\e},\;\, \forall N\ge N_{\e}.
\end{equation*}
\item Let $N\ge N_{\e}$  and let  $h:=\frac{T}{N}$. We denote by $\big(P_{\gk}^N\big)$ the problem $(P_{\gk})$ in which the controls are now restricted to $\U^N$, that is, 
$$ \begin{array}{l} \big(P_{\gk}^N\big)\colon\; \hbox{Minimize}\;g(x(T))\\[1pt] \hspace{1.22cm} \hbox{over}\;(x,u^N)\in AC([0,T];\R^n)\times \mathscr{U}^N\;\hbox{such that for} \;j=1,\dots,N,\\[1pt] \hspace{1.22cm}
 x(\cdot):=x^j(\cdot)\;\hbox{on}\,[(j-1)h, jh]\;\hbox{and}\; x^j \;\hbox{satisfies on}\;  [(j-1)h, jh)\\[3pt]  \hspace{1.2cm} \begin{cases} \dot{x}^j(t)=f(x^j(t),u^N_j)-\gk e^{\gk\psi_{\gk}(x^j(t))} \nabla\psi_{\gk}(x^j(t)), \\[7pt]  x^j((j-1)h)=\begin{sqcases} x_0^{k} & \hbox{if}\;j=1,\\[2pt] x^{j-1}((j-1)h) & \hbox{if}\;j\geq 2. \end{sqcases} \end{cases}
    \end{array}$$
The compactness of $U^N$ yields  that $\big(P_{\gk}^N\big)$ admits an optimal solution. Denote by   $\left(\hat{x}_{\gk}^N,\hat{u}_{\gk}^N\right)$  one of the optimal solutions of $\big(P_{\gk}^N\big)$  and by $\tilde{x}^N_{\gk}$ the solution of $(D)$ corresponding to $\hat{u}_{\gk}^N$. Then, by Proposition \ref{propo2}, we have that 
\begin{equation*}  \label{num3} \|\hat{x}^N_{\gk}-\tilde{x}_{\gk}^N\|_{\infty}\leq \frac{\e}{3L_g},\;\;\; \forall k\ge k_{\e}, \;\, \forall N\ge N_{\e}. \end{equation*} 
\end{enumerate}

The following theorem, Theorem \ref{th1}, is the culmination of all the results of this paper. It basically says that   $\big(P_{\gk}^N\big)$ approximates $(P)$ as $N\f\infty$ and $k\f\infty$. It extends \cite[Theorem 1 $\&$ Remark 4.1]{verachadinum} to our general case, that is, when $C$ is the intersection of a finite number of sublevel sets of smooth functions. Since the statements of our Propositions \ref{propo1}, \ref{propo2} and \ref{propo3} here for the case $r>1$ are, respectively, similar to \cite[Propositions 1, 2 and 3]{verachadinum}, where $r=1$,  then, the proof of Theorem \ref{th1} follows using arguments similar to those used in the proof of  \cite[Theorem 1 $\&$ Remark 4.1]{verachadinum}.

\begin{theorem} \label{th1} For $k\ge k_{\e}$ and $N\ge N_{\e}$, we have  \begin{eqnarray*} &&g(\x(T))\leq g(\tilde{x}^N_{\gk}(T))\leq g(\x(T))+\e,\;\;\hbox{and} \\ && g(\x(T))-\frac{\e}{3}\leq g(\hat{x}^N_{\gk}(T))\leq g(\x(T))+\frac{2\e}{3}.\end{eqnarray*}
Moreover, there exists an optimal solution $(\hat{x},\hat{u})$ of $(P)$ such that, up to a subsequence, both sequences $\hat{x}^N_{\gk}$ and $\tilde{x}^N_{\gk}$ converge uniformly to $\hat{x}$ as $N\f\infty$ and $k\f\infty$.
\end{theorem}

As a consequence of Theorem \ref{th1}, we have the following Algorithm \ref{alg1}  that solves numerically the problem $(P)$.
\begin{algorithm}[htb]
\SetAlgoLined
\KwInput{Positive integer $N$, numbers $\e>0$, $\gamma>0$ and $\delta>0$}\vspace{0.04cm}
\KwOutput{$\blacktriangleright$ Approximating minimum value of $(P)$\newline $\blacktriangleright$ Numerical optimal trajectory of $(P)$} \vspace{0.06cm}
 Initialization\;\vspace{0.04cm}
  $k\gets 1$;\vspace{0.04cm}\\
 \For{$i\gets 0$ \KwTo $1$ \KwBy $1$}{\vspace{0.04cm}
 $\gamma_{i}\gets\gamma+i\sp\delta$;\vspace{0.04cm}\\
 compute $x_0^i$;\vspace{0.04cm}\\
 compute a solution $\hat{x}^N_{\gi}$ of $\big(P_{\gi}^N\big)$ for $\gamma_i$ and $x_0^i$;\vspace{0.04cm}\\
$g_i\gets g\big(\hat{x}^N_{\gi}(T)\big)$;\vspace{0.04cm}\\
}\vspace{0.03cm}
\While{$|g_{k}-g_{k-1}|> \e$}{
  \vspace{0.08cm} $k\gets k+1$;\vspace{0.04cm}\\
  $\gamma_k\gets \gamma_{k-1}+\delta$;\vspace{0.04cm}\\
compute $x_0^k$;\vspace{0.04cm}\\
compute a solution $\hat{x}^N_{\gk}$ of $\big(P_{\gk}^N\big)$ for $\gamma_k$ and $x_0^k$;\vspace{0.04cm}\\
$g_k\gets g\big(\hat{x}^N_{\gk}(T)\big)$;\vspace{0.04cm}\\
 }\vspace{0.03cm}
\Return $g_k$;
\caption{Numerical Method for solving $(P)$}
\label{alg1}
\end{algorithm}

\section{Example} \label{example} To test our algorithm, we provide in this section an example of $(P)$ for which we  separately calculate (i) an exact optimal solution using the Pontryagin-type maximum principle of \cite{verachadijca}, and (ii)  a numerical optimal trajectory using Algorithm \ref{alg1}. Then, we compare our answers.

We consider the following as data for the problem $(P)$:
\begin{enumerate}[$\bullet$]
\item The perturbation mapping $f\colon\R^3\times\R^2\f\R^3$  is defined by $$f((x_1,x_2,x_3),(u,v)):=(x_1-2+u-v,4x_2+x_3+u+v,-x_2+4x_3+u+v).$$
\item The two functions $\psi_1,\;\psi_2\colon \R^3\f\R$ are defined by:
\begin{enumerate}[$\bullet$]
\item $\psi_1(x_1,x_2,x_3):=(x_1-4)^2+x_2^2+x_3^2-25,$ and
\item $\psi_2(x_1,x_2,x_3):=(x_1+4)^2+x_2^2+x_3^2-25.$ 
\end{enumerate}
Hence, the set $C=C_1\cap C_2$ is the intersection of the two solid spheres:
\begin{enumerate}[$\bullet$]
\item $C_1=\{(x_1,x_2,x_3) : (x_1-4)^2+x_2^2+x_3^2\leq 25\}$, and 
\item $C_2=\{(x_1,x_2,x_3) : (x_1+4)^2+x_2^2+x_3^2\leq 25\}$.
\end{enumerate}
Note that $C$ is a {\it nonsmooth} and convex set, see Fig. \ref{Fig1}.
\begin{figure}[ht]
\centering
\includegraphics[scale=0.3]{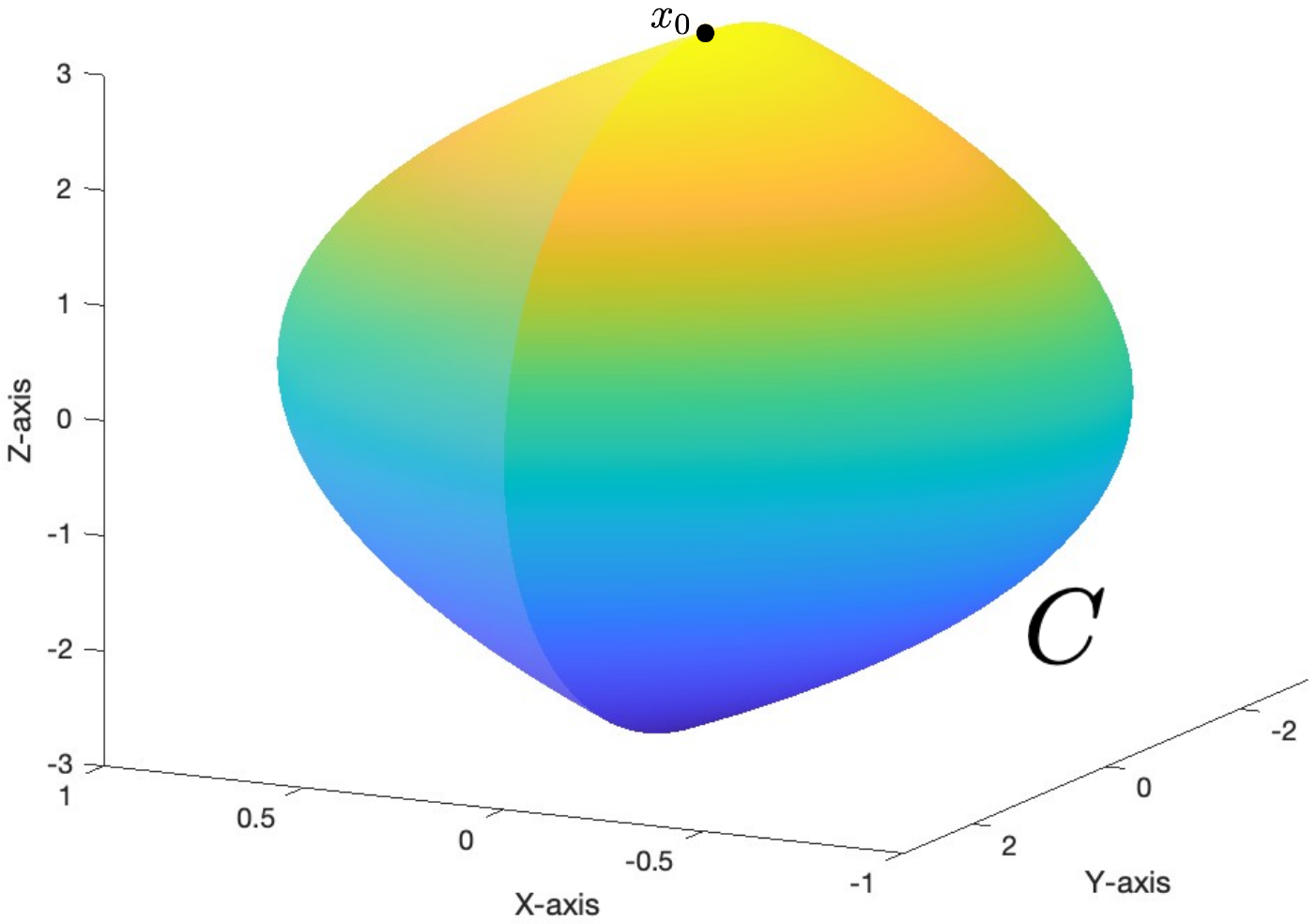}
\caption{\label{Fig1} The sweeping set $C$ of Example \ref{example}}
\end{figure} 
\item The objective function $g\colon\R^3\f\R$  is defined  by $$g(x_1,x_2,x_3):=x_1^2-x_2^2+|x_3|.$$
\item The control set is $U:=[-1,1]\times [-1,1]$, $T:=\tfrac{\pi}{2}$, and the initial point is $x_0:=(0,0,3)$.
\end{enumerate}
One can easily verify that the hypotheses (H1)-(H4), and hence the assumptions (A1)-(A2.2) and (A2.4)-(A6) of  \cite{verachadijca}, are satisfied with $M=35$ and $\eta=2$. Moreover, we have $M_\psi=1$ and $\bar{M}_{\psi}=10$.
\subsection{Exact optimal solution} \label{exactsolution}

In this subsection, we apply the Pontryagin-type maximum principle \cite[Theorem 3.1]{verachadijca} to the problem $(P)$ of our example in order to find one of its optimal solution. Define the curve $$\Gamma:=\{(x_1,x_2,x_3) : x_2^2+x_3^2=9\;\hbox{and}\;x_1=0\}= (\bdry C_1\cap \bdry C_2)\subset\bdry C.$$ 
 Since $x_0\subset \Gamma$ and $g=-x_2^2+|x_3|$ on $\Gamma$ and is strictly greater than $-x^2_2+|x_3|$ elsewhere in $C$, we may seek  for $(P)$  a candidate $(\x,(\u,\vv))$ for  optimality  with  $\x:=(\x_1,\x_2,\x_3)$ belonging to $\Gamma$, if possible,  and hence we have \begin{equation} \label{ex1}\begin{cases} \x_2^2(t)+\x_3^2(t)=9\;\,\hbox{and}\,\;\x_1(t)=0,\;\,\forall t\in [0,\tfrac{\pi}{2}]\,\;\hbox{and}\\[1pt] \x(0)\tran=(0,0,3)\,\;\hbox{and}\;\,\x(\tfrac{\pi}{2})\tran=(0,\pm3,0).\end{cases} \end{equation}
Note that the assumption (A2.3) of \cite{verachadijca} is satisfied on $\Gamma$ for $\beta=(1,1)$.\footnote{For $(x_1,x_2,x_3)\in\Gamma$, we have $\<\nabla\psi_1(x_1,x_2,x_3),\nabla\psi_2(x_1,x_2,x_3)\>=-28<0,$ and hence, the maximum principle of \cite{pinho22} {\it cannot} be applied to this sweeping set $C$.} Then, applying \cite[Theorem 3.1]{verachadijca} to such candidate $(\x,(\u,\vv))$, we obtain the existence of an adjoint vector $p:=(p_1,p_2,p_3)\in BV([0,\frac{\pi}{2}];\R^3)$, two finite signed Radon measures $\nu_1$, $\nu_2$ on $\left[0,\frac{\pi}{2}\right]$, and $\xi_1$, $\xi_2\in L^{\infty}([0,\tfrac{\pi}{2}];\R^{+})$, such that  when incorporating  equations \eqref{ex1}  into \cite[Theorem 3.1$(i)$-$(vi)$]{verachadijca}, we obtain the following: 
\begin{enumerate}[(a)]
\item The admissibility equation holds, that is, for $t\in[0,\tfrac{\pi}{2}]$ a.e.,  \begin{equation*}\begin{cases}0= -2+\u(t)-\vv(t)+8\xi_1(t)-8\xi_2(t),\\ 
\dot{\x}_2(t)= 4\x_2(t)+\x_3(t)+\u(t)+\vv(t)-2\x_2(t)(\xi_1(t)+\xi_2(t)),\\ 
\dot{\x}_3(t)= -\x_2(t)+4\x_3(t)+\u(t)+\vv(t)-2\x_3(t)(\xi_1(t)+\xi_2(t)).\end{cases} \end{equation*} 
\item The adjoint equation is satisfied, that is, for  $t\in [0,\frac{\pi}{2}]$, 
\begin{eqnarray*} dp(t)&=&\begin{pmatrix*}\begin{array}{rrr}
\bp -1\; & 0\;& 0 \\
\bp 0\; & -4\; & 1 \\
\bp 0\; & -1\; & -4 
\end{array}\end{pmatrix*}p(t)\sp dt\sp+\sp (\xi_1(t)+\xi_2(t))\begin{pmatrix} 
2\;\; & 0\;\; & 0  \\
0 \;\; & 2\;\; & 0 \\
0\;\; & 0 \;\; & 2
\end{pmatrix}p(t)\sp dt\sp\\[2pt]&+&\sp \begin{pmatrix} 
-8 \\
2\x_2(t)\\
2\x_3(t)
\end{pmatrix} d\nu_1 + \begin{pmatrix} 
8 \\
2\x_2(t)\\
2\x_3(t)
\end{pmatrix} d\nu_2.\end{eqnarray*}
\item The complementary slackness condition is valid, that is, for $t\in[0,\tfrac{\pi}{2}]$ a.e.,  $$\begin{cases}\xi_1(t)(p_2(t)\x_2(t)+p_3(t)\x_3(t))=0, \\ \xi_2(t)(p_2(t)\x_2(t)+p_3(t)\x_3(t))=0.\end{cases}$$
\item The transversality condition holds, that is, $$p(\tfrac{\pi}{2})\tran\in\{(0,\pm6,\a) : \a\in[-1,1]\}.$$
\item \begin{enumerate}[$\bullet$]
\item $\max\limits_{u\in [-1,1]}\{u\sp(p_1(t)+ p_2(t)+p_3(t))\}$ is attained at $\u(t)$ for $t\in[0,\frac{\pi}{2}]$ a.e.
\item $\max\limits_{v\in [-1,1]}\{v\sp(-p_1(t)+ p_2(t)+p_3(t))\}$ is attained at $\vv(t)$ for $t\in[0,\frac{\pi}{2}]$ a.e.
\end{enumerate}
\end{enumerate}
We temporarily assume that \begin{eqnarray}\label{exam5} \begin{cases} p_1(t)+p_2(t)+p_3(t)\geq 0,\,\;\forall t\in[0,\tfrac{\pi}{2}]\;\hbox{a.e.},\;\hbox{and}\\-p_1(t)+p_2(t)+p_3(t)\leq 0,\,\;\forall t\in[0,\tfrac{\pi}{2}]\;\hbox{a.e.}\end{cases} \end{eqnarray} 
This gives from (e) that $\u(t)=1$ and $\vv(t)=-1$ for $t\in[0,\frac{\pi}{2}]$ a.e. Now solving the differential equations of (a) and using \eqref{ex1}, we obtain that \begin{equation*} \label{lastformu} \xi_1(t)=\xi_2(t)=1\;\,\hbox{and}\;\,\x(t)\tran=(0,3\sin t,3\cos t),\;\;\forall t\in[0,\tfrac{\pi}{2}].\end{equation*}
Hence, from (c) and (d), we deduce that
 \begin{eqnarray}\label{ex4} \begin{cases} p(\tfrac{\pi}{2})\tran\in\{(0,6,\a) : \a\in[-1,1]\},\;\hbox{and}\\ p_2(t)\sin t + p_3(t)\cos t=0,\;\,\forall t\in[0,\tfrac{\pi}{2}]\;\textnormal{a.e.} \end{cases}\end{eqnarray}
Moreover, the adjoint equation (b) simplifies to the following
 \begin{equation}\label{ex3} \begin{cases}d{p}_1(t)= 3p_1(t)\sp dt -8(d\nu_1-d\nu_2),\\
dp_2(t)= p_3(t)\sp dt + 6( d\nu_1+d\nu_2)\sin t,\\
 d{{p}}_3(t)= -p_2(t)\sp dt + 6(d\nu_1+ d\nu_2)\cos t.  \end{cases} \end{equation} 
Using \eqref{ex4} and \eqref{ex3}, a simple calculation gives that 
$$\begin{cases}p(\tfrac{\pi}{2})\tran=(0,6,-1)\;\;\hbox{and}\;\;p(t)\tran=(e^{3t},\cos t,-\sin t)\;\hbox{on}\;[0,\tfrac{\pi}{2}), \\\;d\nu_1=\tfrac{1}{2}\Big(1+\tfrac{e^{\frac{3\pi}{2}}}{8}\Big)\delta_{\left\{\frac{\pi}{2}\right\}}\;\;\hbox{and}\;\;d\nu_2=\tfrac{1}{2}\Big(1-\tfrac{e^{\frac{3\pi}{2}}}{8}\Big)\delta_{\left\{\frac{\pi}{2}\right\}}\;\hbox{on}\;[0,\tfrac{\pi}{2}].\end{cases}$$
where $\delta_{\{a\}}$ denotes the unit measure concentrated on the point $a$. Note that for all $t\in [0,\tfrac{\pi}{2})$, we have $p_1(t)+p_2(t)+p_3(t)\geq 0$ and $-p_1(t)+p_2(t)+p_3(t)\leq 0$, and hence, the temporary assumption \eqref{exam5} is satisfied.

Therefore, the above analysis, realized via  \cite[Theorem 3.1]{verachadijca}, produces  an admissible pair $(\x,(\u,\vv))$, where $$ \x(t)\tran=(0,3\sin t,3\cos t)\;\;\hbox{and}\;\;\u(t)=-\vv(t)=1,\;\;\forall t\in [0,\tfrac{\pi}{2}],$$
which is optimal for $(P)$. This yields that \begin{equation}\label{min(P)}  \min\, (P)=g(0,3,0)=-9.\end{equation}

\subsection{Numerical optimal trajectory} 

The goal of this subsection is to test the effectiveness of Algorithm \ref{alg1}.  Using our algorithm, We  numerically compute  estimates for  both the minimum value and for an optimal trajectory of $\big(P_{\gk}^N\big)$ of  our example. Then, we verify that these computed estimations are in fact good approximations, respectively, for the {\it exact} minimum value and for the {\it exact} optimal trajectory of $(P)$ calculated in Subsection \ref{exactsolution}. This confirms the statement of Theorem \ref{th1}.

\begin{figure}[H]
\centering
\hspace{-0.6cm}
\includegraphics[width=63mm]{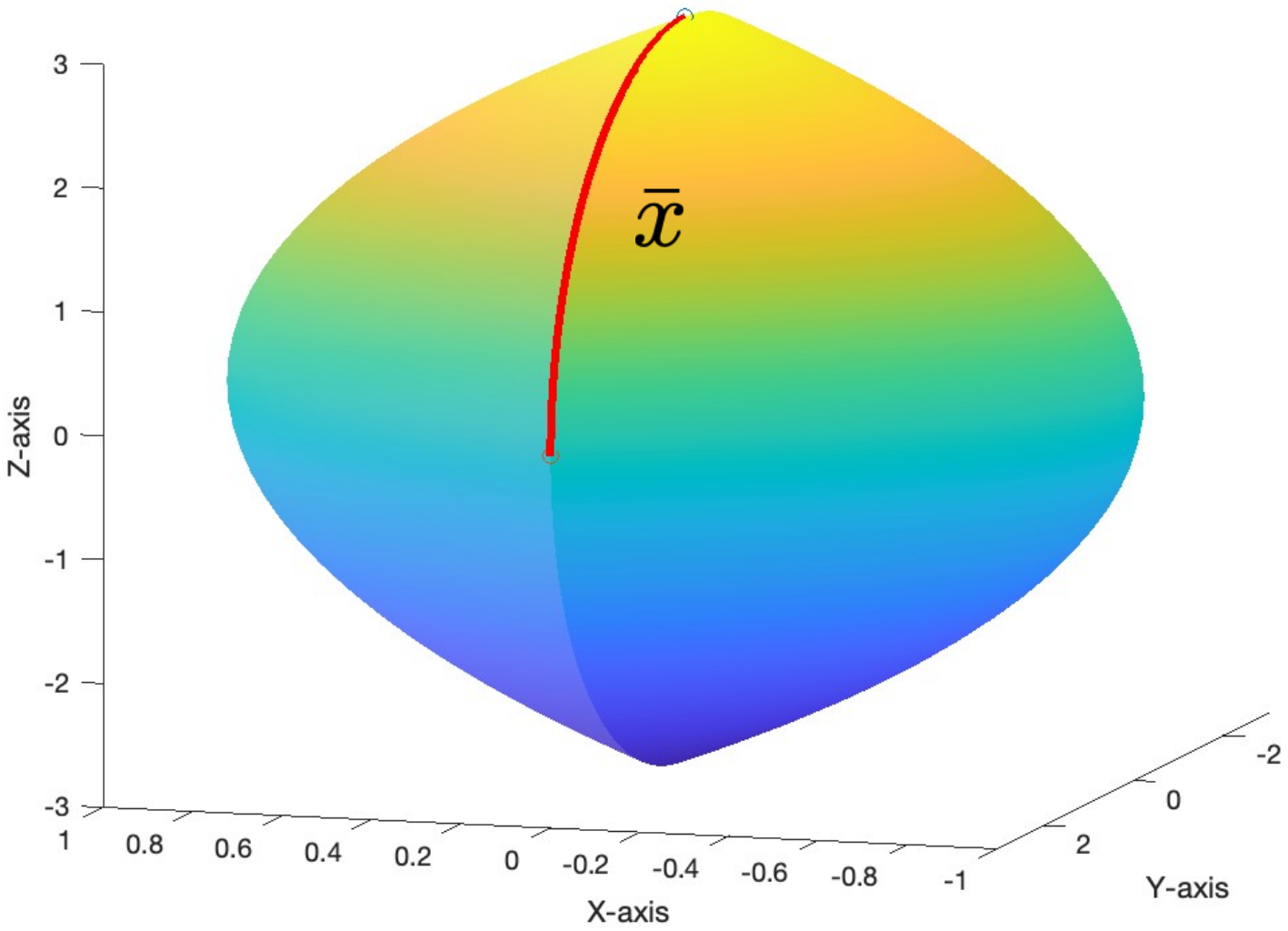}\hspace{-1cm}
\includegraphics[width=69mm]{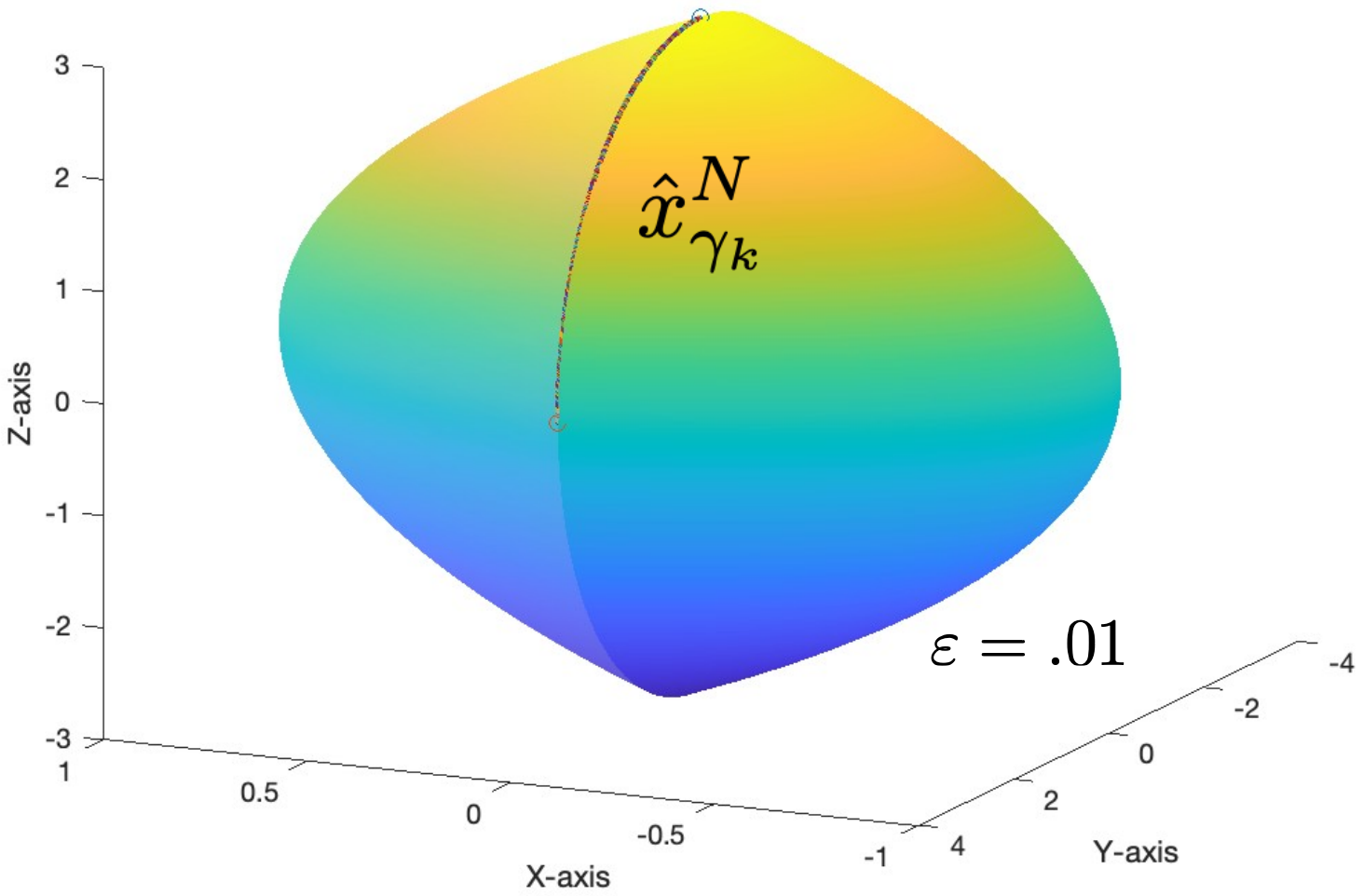}
\includegraphics[width=69mm]{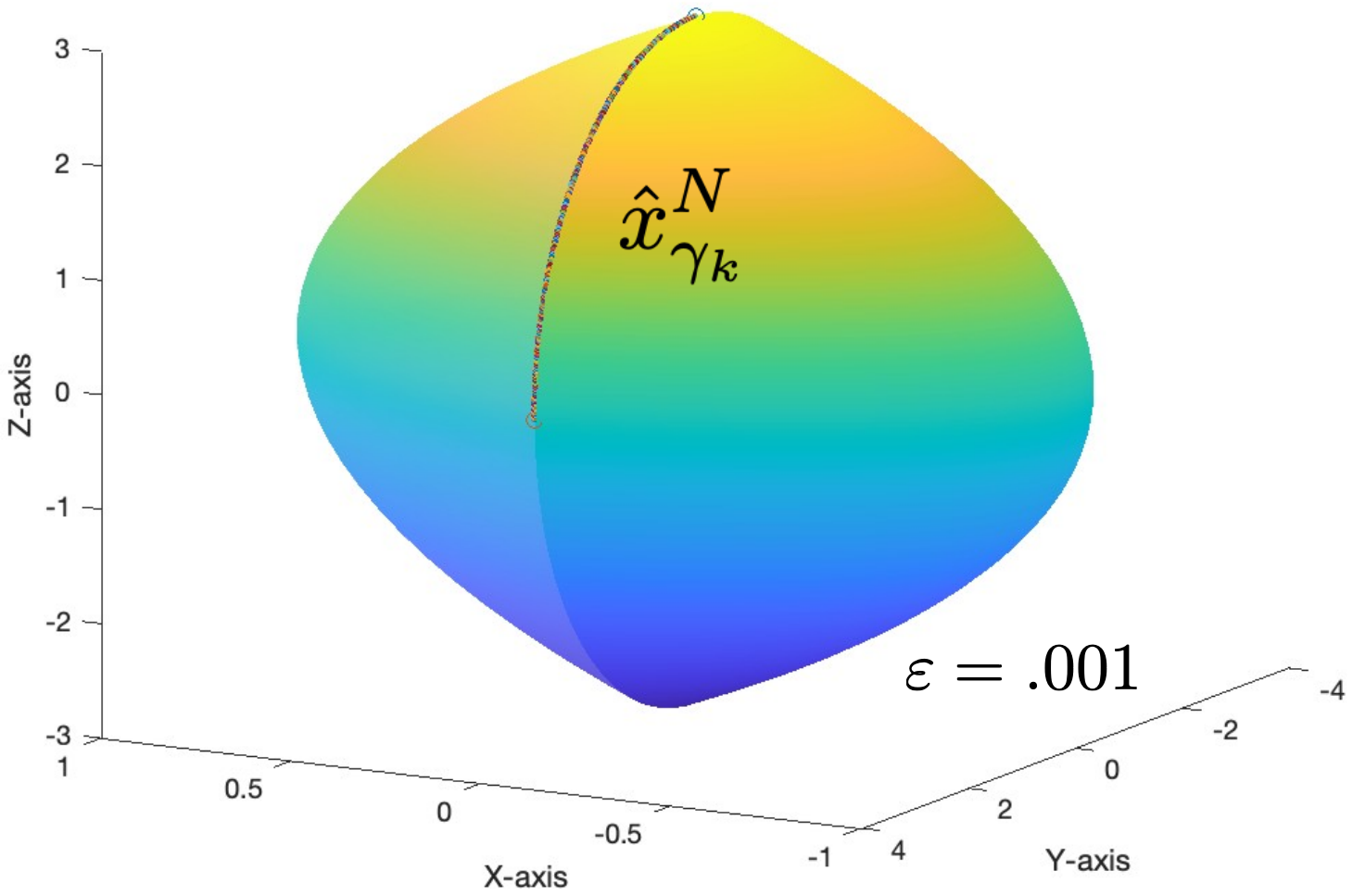}
\caption{\label{Fig2} Exact optimal trajectory $\x$ vs Numerical optimal trajectory $\hat{x}_{\gamma_{k}}^N$}
\vspace{0.2cm}
\end{figure}

We begin by calculating the vector $d_{x_0}$, where $x_0=(0,0,3)$. Clearly we have $\I^0_{x_0}=\{1,2\}$, and hence, $d_{x_0}:=v_1(x_0)+v_2(x_0)$, where $v_1(x_0)$ and $v_2(x_0)$ are the projections of $-\nabla\psi_1(x_0)$ and  $-\nabla\psi_2(x_0)$ to $T_C(x_0)$, respectively. We calculate $N_C(x_0)$ and $T_C(x_0)$, so we find 
\begin{eqnarray*}N_C(x_0)&=&\{\l_1\nabla\psi_1(x_0)+\l_2\nabla\psi_2(x_0) : \l_1\geq 0\;\hbox{and}\;\l_2\geq 0\} \\&=& \{\l_1(-8,,0,6)+\l_2(8,0,6) : \l_1\geq 0\;\hbox{and}\;\l_2\geq 0\},\;\hbox{and}\end{eqnarray*}
$$T_C(x_0)=\{(x_1,x_2,x_3) : [3x_3\leq -4x_1\;\hbox{and}\;x_1\geq0]\;\hbox{or}\;[3x_3\leq 4x_1\;\hbox{and}\;x_1\leq0]\}.$$
Then, the projection of  $-\nabla\psi_1(x_0)=(8,0,-6)$ and $-\nabla\psi_2(x_0)=(-8,0,-6)$ to $T_C(x_0)$ are $\left(\tfrac{144}{25},0,-\tfrac{192}{25}\right)$ and $\left(-\tfrac{144}{25},0,-\tfrac{192}{25}\right)$, respectively. This yields that, $$d_{x_0}= \left(\tfrac{144}{25},0,-\tfrac{192}{25}\right)+\left(-\tfrac{144}{25},0,-\tfrac{192}{25}\right)= \left(0,0,-\tfrac{384}{25}\right).$$
Therefore, \begin{equation}\label{dx0} \frac{d_{x_0}}{\|d_{x_0}\|}= (0,0,-1).\end{equation}
Now, having $M=35$, $\eta=2$ and $M_\psi=1$, we deduce from \eqref{sigmadef}, \eqref{initial} and \eqref{dx0} that for all $k$, 
$$ x_0^k= \left(0,0,3-\tfrac{\ln(2\gk)-\ln 35}{4\gk}\right).$$
We proceed and we write explicitly the problem $\left(P_{\gk}^N\right)$ corresponding to $(P)$.
$$ \begin{array}{l} \big(P_{\gk}^N\big)\colon \hbox{Minimize}\;x^{N}(\frac{\pi}{2})^2-y^N(\tfrac{\pi}{2})^2+|z^N(\tfrac{\pi}{2})|,\vspace{0.1cm}\\ \hspace{1.1cm}\hbox{over}\;((u^1,v^1),\dots,(u^{N},v^{N}))\in ([-1,1]\times[-1,1])^N \vspace{0.1cm}\\\hspace{1.1cm}\hbox{such that for}\; j=1,\dots,N,\;\hbox{we have on}\;\left[(j-1)\frac{\pi}{2N},j\frac{\pi}{2N}\right)\vspace{0.1cm}\\ \hspace{0.95cm} \begin{cases}\dot{x}^j=x^j-2+u^j-v^j-2\xi^1_{\gk}(x^j-4)-2\xi^2_{\gk}(x^j+4),\\[2pt]
\dot{y}^j=z^j+4y^j+u^j+v^j-2\xi^1_{\gk}y^j-2\xi^2_{\gk}y^j,\\[2pt]
\dot{z}^j=-y^j+4z^j+u^j+v^j-2\xi^1_{\gk}z^j-2\xi^2_{\gk}z^j,\\[2.5pt]
\xi^1_{\gk}=\gk e^{\gk((x^j-4)^2+(y^j)^2)+(z^j)^2-25)},\\[2pt]
\xi^2_{\gk}=\gk e^{\gk((x^j+4)^2+(y^j)^2)+(z^j)^2-25)},\\[2pt]
(x^j,y^j,z^j)((j-1)\frac{\pi}{2N})=\begin{sqcases}\left(0,0,3-\tfrac{\ln(2\gk)-\ln 35}{4\gk}\right) & \hbox{if}\;j=1, \\[3pt] (x^{j-1},y^{j-1},z^{j-1})((j-1)\frac{\pi}{2N}) &  \hbox{if}\;j\geq 2. \end{sqcases} 
 \end{cases}
  \end{array}$$ 
We choose $N=20$, $\e=.01$ or $.001$, $\gamma=20$ and $\delta=10$, and we apply Algorithm \ref{alg1}  to numerically compute the minimum value and an approximating optimal trajectory of $(P)$. In order to solve numerically the approximation problem $\big(P_{\gk}^N\big)$, we use MATLAB to implement the Nelder-Mead optimization method  coupled with Runge-Kutta method of fourth order RK4, where the step-size is $.0001$ on each of the 20 intervals. 
\begin{enumerate}[$\bullet$]
\item For $\e=.01$, four iterations of our algorithm reached the desired $\e$ by increasing $\gamma$ to $60$ and the resulting cost is $g_4=g(\hat{x}_{\gamma_4}^N(\tfrac{\pi}{2}))=-8.9319,$ confirming the exact minimum value of $(P)$ found in \eqref{min(P)}. The running time of the algorithm was 75 seconds.\footnote{\label{note1}Machine: MacBook Air, Apple M1 chip, 8GB Unified Memory.}
\item For $\e=.001$, sixteen iterations of our algorithm reached the desired $\e$ by increasing $\gamma$ to $180$ and the resulting cost is $g_{16}=g(\hat{x}_{\gamma_{16}}^N(\tfrac{\pi}{2}))=-8.9693,$ which is now closer to the exact minimum value of $(P)$. The running time of the algorithm was 183 seconds.\footref{note1}
\end{enumerate}
As is easily seen in Fig. \ref{Fig2}  for both cases,    the obtained numerical optimal trajectory $\hat{x}_{\gamma_{k}}^N$ is almost equal to the exact optimal trajectory $\x$ found in Subsection \ref{exactsolution}. This confirms the utility of Theorem \ref{th1}, that is, the convergence, as $N\f\infty$ and $k\f\infty$, of $\hat{x}_{\gamma_{k}}^N$ to an exact optimal trajectory of $(P)$.

\section{Conclusions} In this work, we successfully established  a numerical method to solve optimal control problems involving {\it sweeping processes}, in which the sweeping set is  {\it not} necessarily {\it smooth}, that is,  it is  defined as the intersection of a finite number of sublevel sets of smooth functions. It is worth mentioning that {\it nonsmooth} sweeping sets, including polyhedrals,  are known to occur naturally in applications.

In addition to proving the convergence of our algorithm to an optimal trajectory for the problem, we further confirmed the high effectiveness and efficiency of our numerical method  by providing an example, for which we calculated, on one hand, an exact optimal solution via the maximum principle of \cite{verachadijca},  and, on the other hand, a numerical optimal trajectory via our algorithm. It is remarkable that the approximated trajectory and the exact one turned out to be almost identical, and the error in the objective functions is $.03$,  after running the algorithm for only 183 seconds.

This numerical method  is a  generalization to the nonsmooth setting of the numerical algorithm   in \cite{pinhonum,verachadinum} developed for  smooth sweeping sets.  As opposed to the latter, a number of serious challenges are encountered in this paper.  The nonsmooth property of the maximum function defining our sweeping set posed a major obstacle that prohibited using the technique employed in the smooth setting. To overcome this significant obstacle, we used original and new techniques, including a well constructed smooth approximation of the maximum function defining the sweeping set, and two different, but equivalent, representations of the standard control system that approximates the controlled sweeping process of the original problem. 

Extensions of our numerical method to cover more classes of optimal control problems over sweeping processes, such as  {\it time dependent} nonsmooth sweeping sets, will be the subject of future research.


\end{document}